\RequirePackage{lineno}
 \documentclass[conference,10pt]{IEEEtran}
 \IEEEoverridecommandlockouts
\usepackage{cite,url}
\def\BibTeX{{\rm B\kern-.05em{\sc i\kern-.025em b}\kern-.08em
    T\kern-.1667em\lower.7ex\hbox{E}\kern-.125emX}}
    
\usepackage{threeparttable}
\usepackage{lipsum}
\usepackage{enumitem}
\usepackage[dvips]{color}
\usepackage{lipsum}
\usepackage[ruled,vlined]{algorithm2e}
\usepackage{epsfig,amsmath,amsfonts,cite, enumerate}
\usepackage{booktabs}

\usepackage{multirow}

\DeclareMathAlphabet\mathbfcal{OMS}{cmsy}{b}{n}

\newcommand{\parm}{{\xi}}
\newcommand{\vecpar}{\boldsymbol{\parm}}

\newcommand{\parNum}{d}

\newcommand{\out}{y}

\newcommand{\multiGPC}{\Psi }

\newcommand{\polyInd}{\alpha}
\newcommand{\basisInd}{\boldsymbol{\polyInd}}

\newcommand{\pcOrder}{p}

\newcommand{\yPC}{\sum\limits_{|\basisInd|=0}^{\pcOrder} {c_{\basisInd}  \multiGPC_{\basisInd}  (\vecpar)} }

\newcommand{\parInd}{k}

\newcommand{\mat}[1]{\mathbf{#1}}

\begin{document}
%
\title{ 
Stochastic Collocation with Non-Gaussian Correlated Parameters via a New Quadrature Rule
\thanks{This work was supported by NSF-CCF Award No. 1763699, the UCSB start-up grant and a Samsung Gift Funding.}
}

\author{\IEEEauthorblockN{Chunfeng Cui\IEEEauthorrefmark{1}, Max Gershman\IEEEauthorrefmark{2} and Zheng Zhang\IEEEauthorrefmark{1}}
\IEEEauthorblockA{\IEEEauthorrefmark{1}Department of Electrical and Computer Engineering, University of California Santa Barbara, CA 93106\\
\IEEEauthorrefmark{2}Department of Mathematics/Statistics and Applied Probability,
University of California
Santa Barbara, CA 93106
\\
E-mail:chunfengcui@ucsb.edu; maxwellgershman@umail.ucsb.edu; zhengzhang@ece.ucsb.edu }}

\maketitle

\begin{abstract}
This paper generalizes stochastic collocation methods to handle correlated non-Gaussian random parameters. The key challenge is to perform a multivariate numerical integration in a correlated parameter space when computing the coefficient of each basis function via a projection step. We propose an optimization model and a block coordinate descent solver to compute the required quadrature samples. Our method is verified with a CMOS ring oscillator and an optical ring resonator, showing 3000$\times$ speedup over Monte Carlo.
\end{abstract}


%
\IEEEpeerreviewmaketitle

\section{Introduction}

Stochastic spectral methods are popular techniques to quantify the impact of process variations in nano-scale chip design. Various techniques, such as stochastic Galerkin~\cite{sfem}, stochastic testing~\cite{zzhang:tcad2013}  and stochastic collocation~\cite{col:2005}, have achieved great success in electronic circuits \cite{vrudhula2006hermite, pham2014decoupled,manfredi:tcas2014} and photonics~\cite{weng2015uncertainty}, and have shown significant speedup over Monte Carlo. These techniques approximate a stochastic solution as a linear combination of some basis functions, providing a close-form surrogate model for fast statistical analysis and design automation. 

\begin{figure}[t]
	\centering
\includegraphics[width=3.64in, height=1.2in]{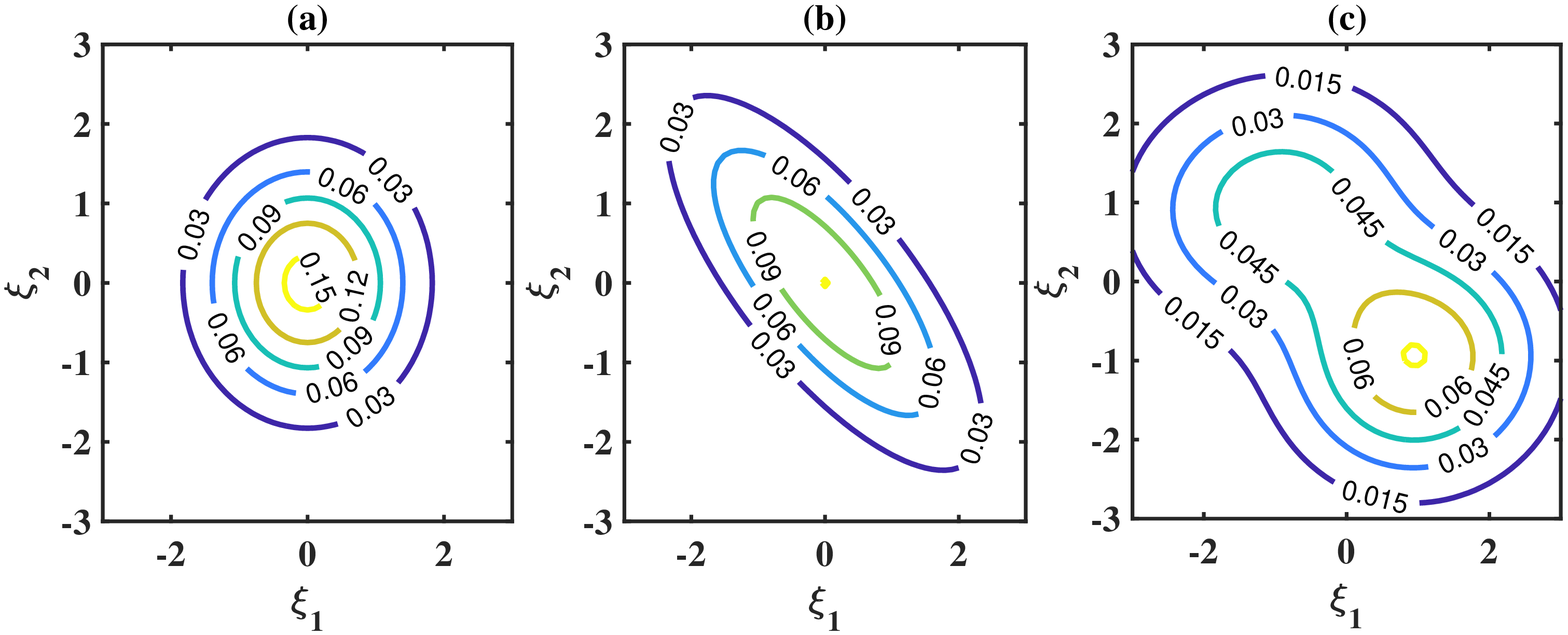}
\caption{Joint density for (a): independent Gaussian, (b): correlated Gaussian, (c): correlated non-Gaussian (e.g., a Gaussian-mixture distribution) cases.}
	\label{fig:gmdistribution}
\end{figure}

Almost all previous stochastic spectral methods assume that the random parameters are mutually independent. This is rarely true in practice. Device geometric or electrical parameters influenced by the same fabrication steps are highly correlated; circuit-level performance parameters used in system-level analysis usually depend on each other.   In this paper, we focus on the non-Gaussian correlated parameters in Fig.~\ref{fig:gmdistribution}~(c). 
Karhunen-Loe\`ve theorem  is error-prone and not scalable.  Preprocessing techniques such as principal component analysis can only handle Gaussian density functions. 


{\bf Our contributions.} We generalize stochastic collocation to non-Gaussian correlated cases by two steps:
\begin{itemize}[leftmargin=*]
\item We propose a new set of basis functions to capture the impact caused by non-Gaussian correlated parameters that cannot be handled by generalized polynomial chaos~\cite{gPC2002}.  
\item Previous integration methods such as sparse grid~\cite{nobile2008sparse} or Gauss quadrature~\cite{Golub:1969} do not work for non-Gaussian correlated cases. Motivated by~\cite{ryu2015extensions,keshavarzzadeh2018numerical}, we propose an optimization solver to calculate the quadrature nodes and weights. We also present a block coordinate descent method to improve the scalability of our solver. 
\end{itemize}
We validate our algorithm by both electronic and photonic ICs, showing $3000\times$ speedup over Monte Carlo.

\section{Review: Stochastic Collocation}


Let $\vecpar=[\parm_1, \cdots, \parm_{\parNum}] \in \mathbb{R}^{\parNum}$ be $d$ random parameters describing process variations. We aim at estimating the uncertainty of a performance metric $y(\vecpar)$ (e.g., chip frequency or power).
Stochastic spectral methods approximate the solution by
\begin{equation}
\label{eq:ygpc}
\out (\vecpar) \approx \yPC, \; {\rm with}\; \mathbb{E}\left[{\multiGPC}_{\basisInd} \left( \vecpar \right)\multiGPC_{\boldsymbol{\beta }}\left( \vecpar \right)\right ]=\delta_{\basisInd, \boldsymbol{\beta }}.
 \end{equation}
Here $\mathbb{E}$ denotes expectation, $\delta$ denotes a Delta function, the basis functions $\{{\multiGPC}_{\basisInd} \left(\vecpar\right)\}$ are  orthonormal polynomials, $\basisInd=[\alpha_1,\cdots, \alpha_{\parNum}] \in \mathbb{N}^{\parNum}$ indicates the highest polynomial order of each parameter in the corresponding basis. The total polynomial order $|\basisInd|=\alpha_1+\ldots+\alpha_d$ is bounded by $p$, and thus the total number of basis functions is $
N=(p+d)!/(p!d!)$.

Projection-based stochastic collocation methods compute the coefficient $c_{\basisInd}$ via a numerical integration. If one has $M$ quadrature nodes  $\{\vecpar_{\parInd}\}_{k=1}^M$ and weights $\{w_{\parInd}\}_{k=1}^M$, then
\begin{equation}
\label{eq:yTP}
 c_{\basisInd}=\mathbb{E} \left[ \out (\vecpar)  {\multiGPC}_{\basisInd} (\vecpar)\right]\approx \sum\limits_{\parInd=1}^M {\out (\vecpar_{\parInd})  {\multiGPC}_{\basisInd}  (\vecpar_{\parInd})  w_{\parInd}  }.
\end{equation}
If $\vecpar$ are mutually independent, then ${\multiGPC}_{\basisInd}  (\vecpar_{\parInd})$ may be chosen as the generalized polynomial chaos~\cite{gPC2002}, and the quadrature nodes and weights can be calculated via sparse grid~\cite{nobile2008sparse} and Gauss quadrature~\cite{Golub:1969}. However, how to choose the basis functions and quadrature rule is an open question for non-Gaussian correlated cases. Soize suggested a modification of generalized polynomial chaos~\cite{Soize2004physical}, but the resulting basis functions are non-smooth and unstable~\cite{Cui2018}.

 \section{Our basis functions}


We adopt the Gram-Schmidt approach to calculate the basis function recursively.
 Gram-Schmidt was originally used for vector orthogonalization in the Euclidean space, and the key difference here is to replace the vector inner product with a functional expectation. Specifically, we first reorder the monomials $\vecpar^{\basisInd}=\xi_1^{\alpha_1}\ldots \xi_d^{\alpha_d}$ in the graded lexicographic order, and denote them as $\{p_j(\vecpar)\}_{j=1}^N$.
 Then we set $\multiGPC_1(\vecpar)   = 1$ and calculate a set of orthonormal polynomials $\{\multiGPC_{j}(\vecpar)\}_{j=1}^N$ in the correlated parameter space recursively:
 \begin{align} &\hat{\multiGPC}_j(\vecpar) = p_j(\vecpar)-\sum_{i=1}^{j-1} \mathbb{E}[ p_j(\vecpar)\multiGPC_i(\vecpar)] \multiGPC_i(\vecpar),
 \\
 &\multiGPC_j(\vecpar)  = \frac{\hat{\multiGPC}_j(\vecpar)}{\sqrt{\mathbb{E}[\hat{\multiGPC}^2_j(\vecpar)]}},\ j=2,\ldots,N.
 \end{align}
The most time-consuming step is to compute the expectations. We adopt the  functional tensor train approach developed in~\cite{Cui2018} to speed up this computation.



\section{An Optimization-Based Quadrature Rule}
Having chosen the basis functions, we still need to determine the number and values of the quadrature nodes and weights in order to calculate the coefficient $c_{\basisInd}$ by (\ref{eq:yTP}). Our proposed method is summarized in Algorithm~\ref{alg:updatenodes}, and we explain the key ideas as follows.
 

 \subsection{An Optimization-Based Quadrature Rule}
 \label{subsec:BCD}

Motivated by \cite{ryu2015extensions, keshavarzzadeh2018numerical}, we set up an optimization model to decide a proper quadrature rule. Our method differs from~\cite{ryu2015extensions} because the latter optimizes quadrature weights only. Our method differs from~\cite{keshavarzzadeh2018numerical} in the following sense: (1) we focus on non-Gaussian correlated uncertainty analysis; (2) we handle the  nonnegative constraint of $\mat{w}$ and the nonlinear function of $\vecpar$ separately via a novel block coordinate descent framework.  

Suppose that $\out (\vecpar)$ can be well approximated by an order-$p$ polynomial function, then the product term $ \out (\vecpar)  {\multiGPC}_{\basisInd} (\vecpar)$ can be well approximated  by order-$2p$ polynomials. As a result, $\mathbb{E} \left[ \out (\vecpar)  {\multiGPC}_{\basisInd} (\vecpar)\right]$ can be accurately computed if we have a quadrature rule that can accurately estimate the integration of {\it every} basis function bounded by order $2p$:
 \begin{align}
 \label{equ:nmint1}
 \mathbb{E}[\multiGPC_{j}(\vecpar)]=&\delta_{1j} \approx\sum_{k=1}^{M}\multiGPC_j(\vecpar_k)  w_k,\ \forall\, j=1,\ldots,N_{2p},
 \end{align}
with $N_{2p}=\binom{d+2p}{d}$,  $\delta_{1j}=1$ if $j=1$ and $\delta_{1j}=0$ otherwise.
This formulation can be rewritten as a nonlinear least-square 
 \begin{equation}\label{equ:NLS}
\min_{\bar{\vecpar},\mat{w}}\quad \|\mat{\Phi}(\bar{\vecpar})\mat{w}-\mat{e_1}\|^2,
\end{equation}
where $(\mat{\Phi}(\bar{\vecpar}))_{jk}=\multiGPC_j(\vecpar_k)$, $\bar{\vecpar}=[\vecpar_1;\ldots;\vecpar_M]\in\mathbb{R}^{Md}$, $\mat{w}=[w_1,\ldots,w_M]^T\in\mathbb{R}^{M}$ and $\mat{e}_1=[1,0,\ldots,0]^T\in\mathbb{R}^{N_{2p}}$.

\subsection{A Block Coordinate Solver for \eqref{equ:NLS}}

The number of unknowns in \eqref{equ:NLS} is $M(d+1)$, which becomes large as $d$ increases. To improve the scalability, we solve \eqref{equ:NLS}  by a block coordinate descent method. The idea is to update the variables block-by-block: at the $t$-th  iteration, given $\bar{\vecpar}_t$ and $\mat{w}_t$, we firstly fix $\bar{\vecpar}^t$ and solve the $\mat{w}$-subproblem to update $\mat{w}^{t+1}$, then fix $\mat{w}^{t+1}$ and solve the $\vecpar$-subproblem to get $\bar{\vecpar}^{t+1}$.

 \textbf{$\mat{w}$-subproblem.} Suppose $\bar{\vecpar}^t=[\vecpar_1^t;\ldots;\vecpar_M^t]$ is fixed, then (\ref{equ:NLS}) reduces to a convex linear least-square problem
  \begin{equation}\label{equ:NLSw}
\mat{w}^{t+1}=\arg\min_{\mat{w}\ge0}\quad \|\mat{\Phi}(\bar{\vecpar}^t) \mat{w}-\mat{e}_1\|^2.
 \end{equation}
Here, we require the quadrature weights  to be nonnegative.

 \textbf{$\vecpar$-subproblem.} When $\mat{w}^{t+1}$ is fixed, we apply the Gaussian Newton method to  the $\vecpar$-subproblem 
 \begin{equation*}
 \vecpar_k^{t+1}=\vecpar_k^t+\mat{d}_k^t, \text{ where } \{\mat{d}_k^t\}=\arg\min_{\{\mat{d}_k\}} \ \|\sum_{k=1}^M\mat{G}_k^t\mat{d}_k + \mat{r}^t\|^2.
 \end{equation*}
 Here,  $\mat{r}^t = \mat{\Phi}(\bar{\vecpar}^t)\mat{w}^{t+1}-\mat{e}_1 \in \mathbb{R}^{N_{2p}}$ denotes the residual, $\mat{G}_k^t\in \mathbb{R}^{N_{2p}\times d}$ is the Jacobian matrix of $\mat{r}^t$ with respect to $\vecpar_k^t$. 
 

 \begin{algorithm}[t]
\label{alg:updatenodes}
\caption{Extensions of stochastic collocation method to non-Gaussian correlated variations}
      \SetKwInput{Input}{Input}
      \SetKwInput{Output}{Output}
 \begin{itemize}
  \item[Step 1] Initialize the quadrature nodes and weight according to Section \ref{section:initialNodes}.
 \item[Step 2] \textbf{Increase phase.} Update the quadrature nodes \& weights by solving~\eqref{equ:NLS}. If the optimization fails to converge, increase the node number and go back to Step~1.
 \item[Step 3] \textbf{Decrease phase.}  Decrease the node number, and update them by solving~\eqref{equ:NLS} again. Repeat Step~3 until no points can be deleted. Return the optimal nodes and weights. 
 \item[Step 4] Call a simulator to compute $\{\out(\vecpar_k)\}_{k=1}^M$. Then compute the coefficients $\{c_{\basisInd}\}$ for all ${|\basisInd|\le p}$ via (\ref{eq:yTP}).
 \end{itemize}
 \Output{The coefficients $\{c_{\basisInd}\}$ in (\ref{eq:ygpc}).}
\end{algorithm}

 \subsection{Implementation Details}
 \label{section:initialNodes}
A good initial guess for the quadrature nodes is important to ensure the success of our nonlinear least-square solver. Therefore, we first generate some candidate nodes via a Monte Carlo method, and then cluster them   via a
complete-linkage clustering method~\cite{defays1977efficient}.

In general, we do not know the optimal number of quadrature nodes \textsl{a priori}. Our algorithm consists of two phases: firstly we increase the number of quadrature nodes until the  condition (\ref{equ:nmint1}) holds with high accuracy. Then, we decrease the number of nodes by deleting the node with the least weight and refine them by solving~\eqref{equ:NLS}, until the number of nodes is too small to achieve a required integration accuracy. 


\section{Numerical Results}

\subsection{Three-Stage CMOS Ring Oscillator}
\begin{figure}[t]
	\centering
		\includegraphics[width=2.4in]{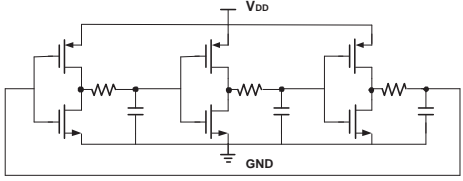}
\caption{Schematic of a 3-stage CMOS ring oscillator. }
	\label{fig:ring}
\end{figure}

\begin{figure}[t]
	\centering
	\includegraphics[width=85mm, height=1.5in]{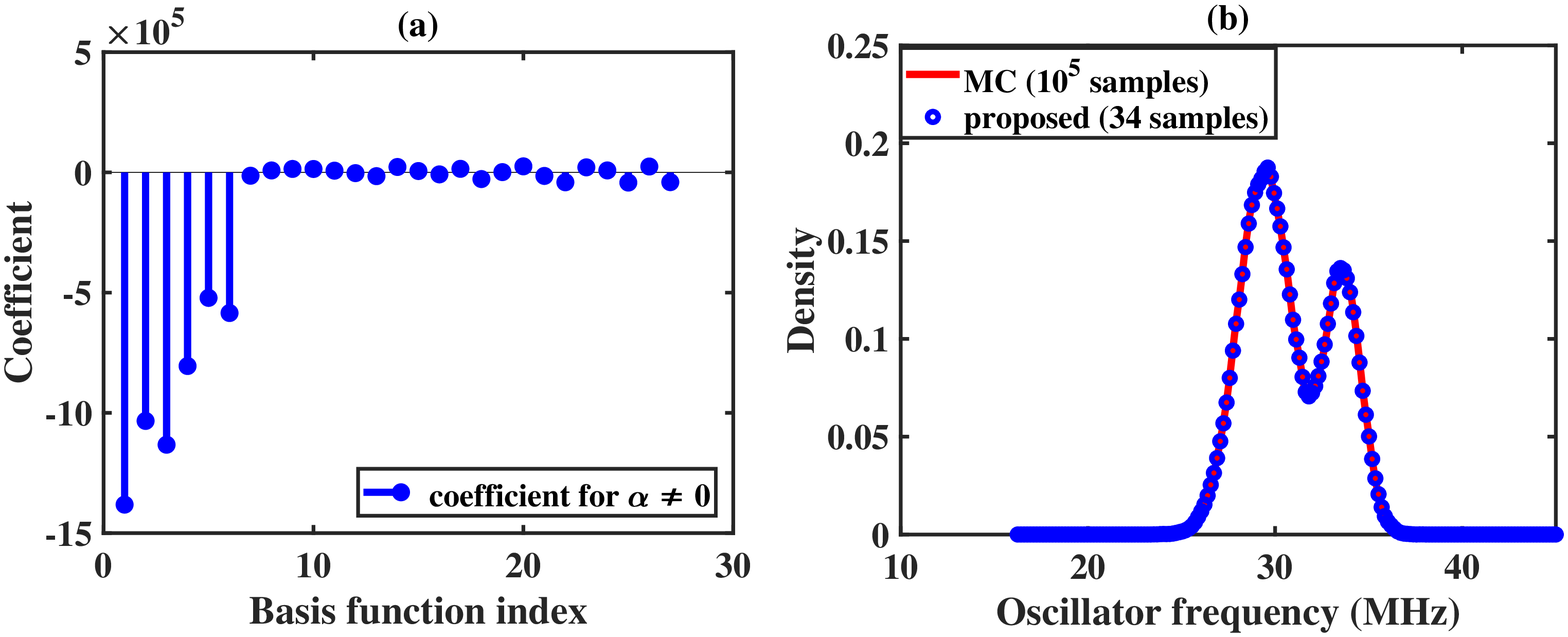}
\caption{Numerical results of the CMOS ring oscillator. (a) obtained coefficients/weights of our basis functions; (b) probability density functions of the oscillator frequency obtained by our proposed method and Monte Carlo (MC).}
	\label{fig:ring_results}
\end{figure}

We first use our method to simulate the 3-stage CMOS ring oscillator in Fig.~\ref{fig:ring}. This oscillator has a Gaussian mixture model describing the correlated non-Gaussian threshold voltages of 6 transistors. We aim to obtain a 2nd-order expansion for its frequency by calling a periodic steady-state simulator repeatedly. The obtained results in Fig.~\ref{fig:ring_results}  
shows the obtained coefficients for all basis functions. The obtained density function using only $34$ quadrature samples is almost identical with that from $10^5$ Monte Carlo simulations.


\subsection{Coupled Ring Resonator Optical Filter}
\label{subsec:ringresonator}

\begin{figure}[t] 
\centering
\includegraphics[width=1.6in, height=1in]{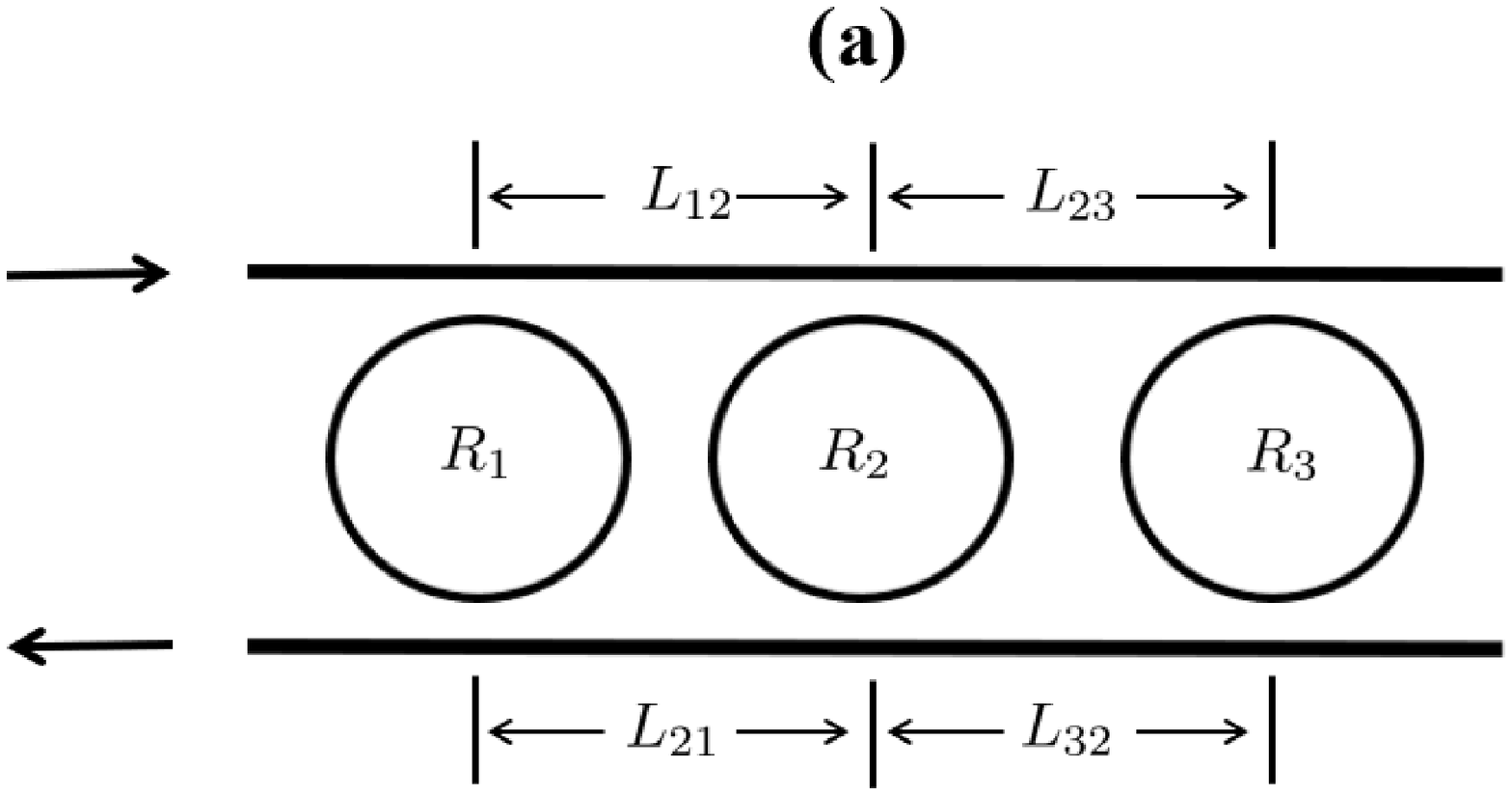}
\includegraphics[width=1.8in, height=1.2in]{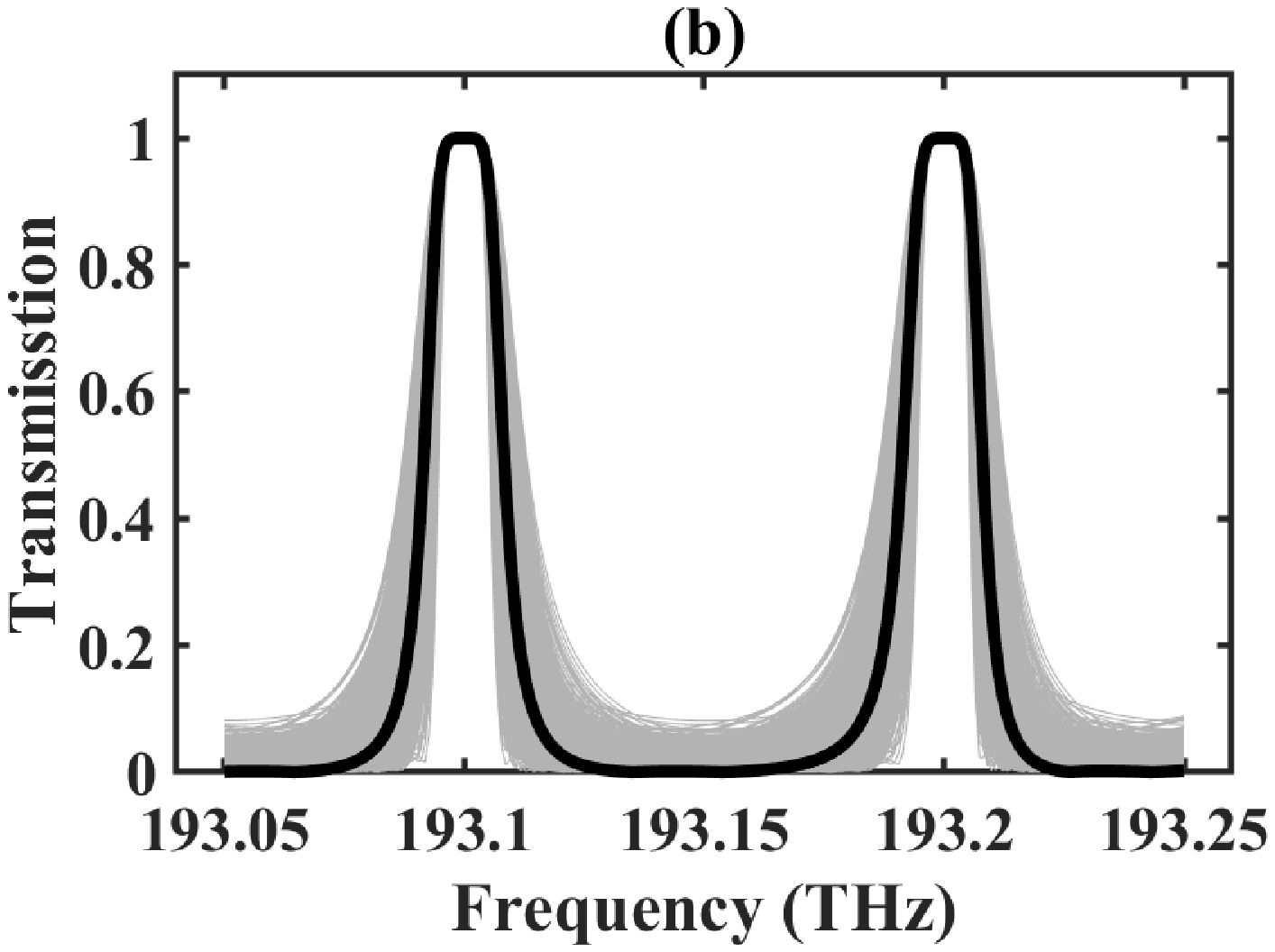}
\caption{(a) Schematic of a 3-stage parallel-coupled ring resonator optical filter. 
(b) The black line shows the nominal transmission function, and the thin grey lines show the effect of fabrication uncertainties on the waveguide lengths.
}
\label{fig:coupled_ring_resonator}
\end{figure}

We further consider the filter designed with bus-coupled micro-ring resonators\footnote{The details of this benchmark can be found at \url{https://kb.lumerical.com/en/pic_circuits_coupled_ring_resonator_filters.html}} shown in Fig.~\ref{fig:coupled_ring_resonator}~(a). Coupled ring resonator are widely used for wavelength filtering and modulation in photonic integrated circuits. Here we consider a filter with $3$ stages of ring resonators, and we use a Gaussian mixture model to describe the correlated non-Gaussian uncertainties in waveguide lengths $L_{12}$, $L_{21}$, $L_{23}$ and $L_{32}$. 

\begin{figure}[t] 
	\centering 
	\includegraphics[width=90mm,height=1.5in]{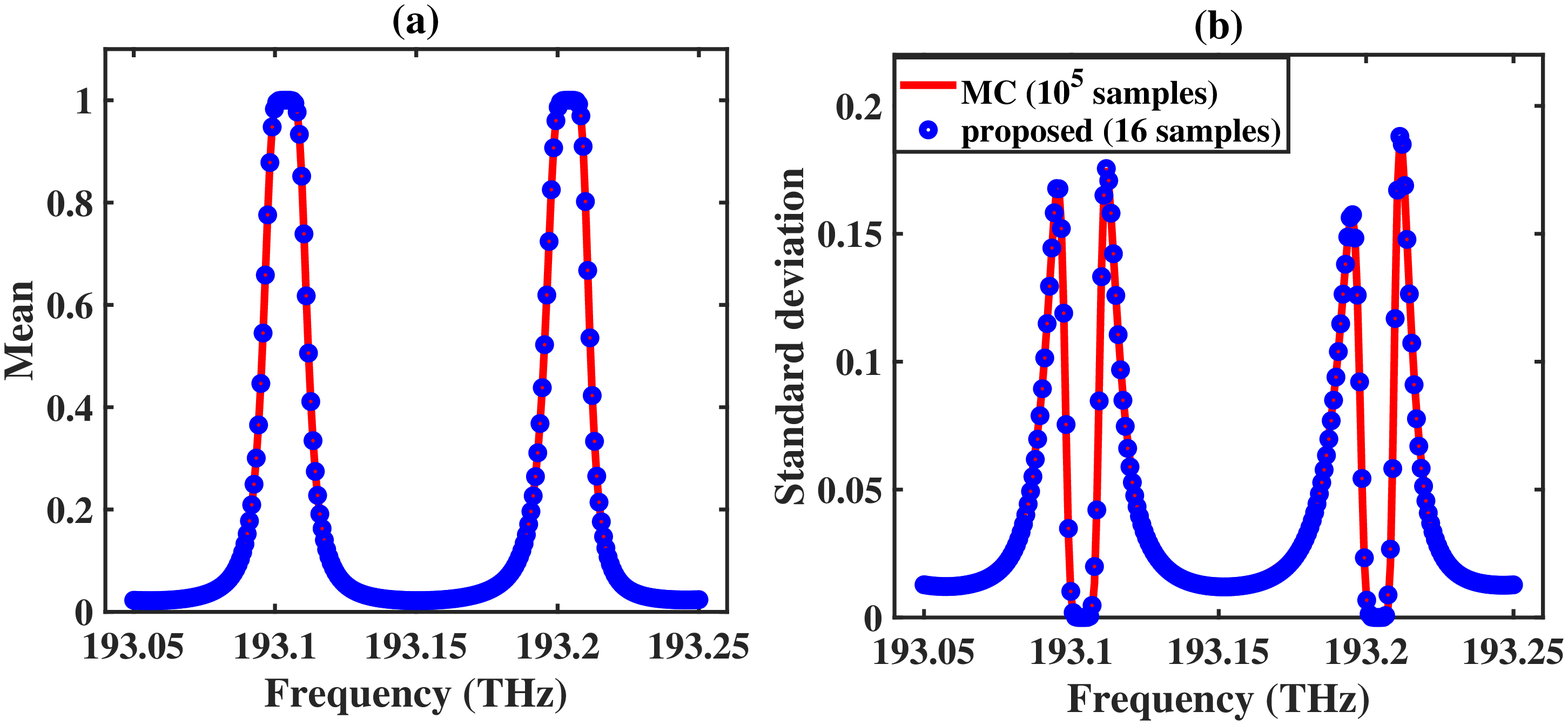}
\caption{Simulation results of the optical filter. (a) obtained mean value of the power transmission rate; (b) standard deviation of the transmission rate.} \label{fig:res_ring2}
\end{figure}

A 2nd-order expansion is built to approximate the frequency-dependent power transmission function: $\out(f,\vecpar)=\sum_{|\basisInd|=0}^p c_{\basisInd}(f) \multiGPC_{\basisInd}(\vecpar)$.
The computed mean value and standard derivation are shown in Fig.~\ref{fig:res_ring2}. Our method only uses 16 quadrature samples for simulation, and it is able to achieve the similar level of accuracy compared with Monte Carlo method using $10^5$ simulation samples.

\section{Conclusion}
This paper has proposed a stochastic collocation approach to solve the challenging non-Gaussian correlated uncertainty quantification problems.
We have proposed an optimization method to calculate the quadrature rule used in the projection step. Our method has achieved 3000$\times$ speedup than Monte Carlo on a CMOS ring oscillator and an optical resonator.




\bibliographystyle{IEEEtran}
\bibliography{EPEPS}

\begin{thebibliography}{10}
\providecommand{\url}[1]{#1}
\csname url@samestyle\endcsname
\providecommand{\newblock}{\relax}
\providecommand{\bibinfo}[2]{#2}
\providecommand{\BIBentrySTDinterwordspacing}{\spaceskip=0pt\relax}
\providecommand{\BIBentryALTinterwordstretchfactor}{4}
\providecommand{\BIBentryALTinterwordspacing}{\spaceskip=\fontdimen2\font plus
\BIBentryALTinterwordstretchfactor\fontdimen3\font minus
  \fontdimen4\font\relax}
\providecommand{\BIBforeignlanguage}[2]{{%
\expandafter\ifx\csname l@#1\endcsname\relax
\typeout{** WARNING: IEEEtran.bst: No hyphenation pattern has been}%
\typeout{** loaded for the language `#1'. Using the pattern for}%
\typeout{** the default language instead.}%
\else
\language=\csname l@#1\endcsname
\fi
#2}}
\providecommand{\BIBdecl}{\relax}
\BIBdecl

\bibitem{sfem}
R.~Ghanem and P.~Spanos, \emph{Stochastic finite elements: a spectral
  approach}.\hskip 1em plus 0.5em minus 0.4em\relax Springer-Verlag, 1991.

\bibitem{zzhang:tcad2013}
Z.~Zhang, T.~A. El-Moselhy, I.~A.~M. Elfadel, and L.~Daniel, ``Stochastic
  testing method for transistor-level uncertainty quantification based on
  generalized polynomial chaos,'' \emph{IEEE Trans. Computer-Aided Design
  Integr. Circuits Syst.}, vol.~32, no.~10, Oct. 2013.

\bibitem{col:2005}
D.~Xiu and J.~S. Hesthaven, ``High-order collocation methods for differential
  equations with random inputs,'' \emph{SIAM J. Sci. Comp.}, vol.~27, no.~3,
  pp. 1118--1139, Mar 2005.

\bibitem{vrudhula2006hermite}
S.~Vrudhula, J.~M. Wang, and P.~Ghanta, ``Hermite polynomial based interconnect
  analysis in the presence of process variations,'' \emph{IEEE Trans. CAD of
  Integr. Circuits Syst.}, vol.~25, no.~10, pp. 2001--2011, 2006.

\bibitem{pham2014decoupled}
T.-A. Pham, E.~Gad, M.~S. Nakhla, and R.~Achar, ``Decoupled polynomial chaos
  and its applications to statistical analysis of high-speed interconnects,''
  \emph{IEEE Transactions on Components, Packaging and Manufacturing
  Technology}, vol.~4, no.~10, pp. 1634--1647, 2014.

\bibitem{manfredi:tcas2014}
P.~Manfredi, D.~V. Ginste, D.~D. Zutter, and F.~Canavero, ``Stochastic modeling
  of nonlinear circuits via {SPICE}-compatible spectral equivalents,''
  \emph{IEEE Trans. Circuits Syst. I: Reg. Papers}, vol.~61, no.~7, pp.
  2057--2065, July 2014.

\bibitem{weng2015uncertainty}
T.-W. Weng, Z.~Zhang, Z.~Su, Y.~Marzouk, A.~Melloni, and L.~Daniel,
  ``Uncertainty quantification of silicon photonic devices with correlated and
  non-{Gaussian} random parameters,'' \emph{Optics Express}, vol.~23, no.~4,
  pp. 4242--4254, 2015.

\bibitem{gPC2002}
D.~Xiu and G.~E. Karniadakis, ``The {Wiener-Askey} polynomial chaos for
  stochastic differential equations,'' \emph{SIAM J. Sci. Comp.}, vol.~24,
  no.~2, pp. 619--644, Feb 2002.

\bibitem{nobile2008sparse}
F.~Nobile, R.~Tempone, and C.~G. Webster, ``A sparse grid stochastic
  collocation method for partial differential equations with random input
  data,'' \emph{SIAM J. Numerical Analysis}, vol.~46, no.~5, pp. 2309--2345,
  2008.

\bibitem{Golub:1969}
G.~H. Golub and J.~H. Welsch, ``Calculation of {Gauss} quadrature rules,''
  \emph{Math. Comp.}, vol.~23, pp. 221--230, 1969.

\bibitem{ryu2015extensions}
E.~K. Ryu and S.~P. Boyd, ``Extensions of {Gauss} quadrature via linear
  programming,'' \emph{Foundations of Computational Mathematics}, vol.~15,
  no.~4, pp. 953--971, 2015.

\bibitem{keshavarzzadeh2018numerical}
V.~Keshavarzzadeh, R.~M. Kirby, and A.~Narayan, ``Numerical integration in
  multiple dimensions with designed quadrature,'' \emph{arXiv preprint
  arXiv:1804.06501}, 2018.

\bibitem{Soize2004physical}
C.~Soize and R.~Ghanem, ``Physical systems with random uncertainties: chaos
  representations with arbitrary probability measure,'' \emph{SIAM Journal on
  Scientific Computing}, vol.~26, no.~2, pp. 395--410, 2004.

\bibitem{Cui2018}
C.~Cui and Z.~Zhang, ``Uncertainty quantification of electronic and photonic
  {ICs} with non-{Gaussian} correlated process variations,'' \emph{arXiv
  preprint arXiv:1807.01778}, 2018.

\bibitem{defays1977efficient}
D.~Defays, ``An efficient algorithm for a complete link method,'' \emph{The
  Computer Journal}, vol.~20, no.~4, pp. 364--366, 1977.

\end{thebibliography}

\end{document}